\documentclass[leqno,10pt]{article}
\usepackage{amssymb,amsmath,amsfonts,amsthm,latexsym,ifthen,makeidx,bezier}
\usepackage{graphicx}
\usepackage{undertilde}
\usepackage{amsxtra}
\usepackage[a4paper]{geometry}
\usepackage{url}
\usepackage[english]{babel} 
\usepackage[tight,centredisplay]{diagrams}

\newcommand{\kla}[1]{ {\langle #1 \rangle} }
\newcommand{\anf}[1]{{\text{``}{#1}\text{''}}}

\newcommand{\st}{\;|\;}

\newcommand{\dom}{ {\rm dom} }

\newcommand{\sub}{\subseteq}

\newfont{\ssi}{cmssi12 at 12pt}

\newcommand{\On}{ {\rm On} }

\newcommand{\leer}{\emptyset}
\newcommand{\ohne}{\setminus}

\newenvironment{ea*}{\begin{eqnarray*}}{\end{eqnarray*}}

\setbox0=\hbox{$\longrightarrow$}
\newcommand{\To}{\longrightarrow}

\newcommand{\vc}{{\vec{c}}}

\newcommand{\vp}{{\vec{p}}}

\newcommand{\valpha}{{\vec{\alpha}}}

\newcommand{\seq}[2]{{\langle#1\;|\;}\linebreak[0]{#2\rangle}}

\renewcommand{\phi}{\varphi}

\newcommand{\ZFC}{\ensuremath{\mathsf{ZFC}}}
\newcommand{\ZF}{\ensuremath{\mathsf{ZF}}}

\newcommand{\V}{\ensuremath{\mathrm{V}}}

\newcommand{\forces}{\Vdash}

\def\<#1>{\langle#1\rangle}

\renewcommand{\P}{{\mathord{\mathbb P}}}

\newcommand{\of}{\subseteq}

\newcommand{\MP}{\ensuremath{\mathsf{MP}}}

\newcommand{\ColNothing}{\mathrm{Col}}
\newcommand{\Col}[1]{\ColNothing(#1)}

\newcommand{\MPColNothing}[1]{\MP_{\Col{\dot{\kappa}}}}

\newcommand{\OD}{\ensuremath{\mathsf{OD}}}
\newcommand{\HOD}{\ensuremath{\mathsf{HOD}}}

\newcommand{\reals}{{\mathord{\mathbb{R}}}}

\newcommand{\isomorphism}{\stackrel{\sim}{\longleftrightarrow}}

\newtheorem{thm}{Theorem}[section]
\newtheorem{cor}[thm]{Corollary}
\newtheorem{lem}[thm]{Lemma}
\newtheorem{obs}[thm]{Observation}

\theoremstyle{definition}
\newtheorem{defn}[thm]{Definition}
\newtheorem{question}[thm]{Question}

\theoremstyle{remark}

\usepackage{verbatim}
\usepackage{tikz} 
\usepackage[hidelinks]{hyperref}
\usetikzlibrary{arrows}
\usepackage{graphicx}
\newcommand{\type}{\mathrm{tp}}
\newcommand{\Ehrenfeucht}{\mathsf{EL}}

\newcommand{\ExternalAlgebraicEhrenfeucht}{\mathsf{AEL}_{\text{ext}}}
\newcommand{\InternalAlgebraicEhrenfeucht}{\mathsf{AEL}_{\text{int}}}
\newcommand{\InternalOrdinalAlgebraic}{\mathsf{OA}_{\text{int}}}
\newcommand{\ExternalOrdinalAlgebraic}{\mathsf{OA}_{\text{ext}}}
\newcommand{\OrdinalAlgebraic}{\mathsf{OA}}
\newcommand{\OA}{\OrdinalAlgebraic}
\newcommand{\HOA}{\mathsf{HOA}}
\newcommand{\ExternalAlgebraic}{\mathsf{A}_{\text{ext}}}
\newcommand{\InternalAlgebraic}{\mathsf{A}_{\text{int}}}
\newcommand{\Definable}{\mathsf{D}}
\newcommand{\PA}{\ensuremath{\mathsf{PA}}}
\newcommand{\image}{\mathbin{\hbox{\tt\char'42}}}
\newcommand{\LM}{\ensuremath{\mathsf{LM}}}
\newcommand{\Sacks}{\mathbb{S}}
\newcommand{\cdegree}[1]{[\!\!|#1|\!\!]}

\begin{document}

\title{Ehrenfeucht's lemma in set theory}

\author{Gunter Fuchs, Victoria Gitman and
        Joel David Hamkins\thanks{The research of the third author has been supported
        in part by NSF grant DMS-0800762, PSC-CUNY grant 64732-00-42, CUNY Collaborative
        Incentive Award 80209-06 20 and Simons Foundation grant 209252.\newline Commentary concerning this article can be made at \href{http://jdh.hamkins.org/ehrenfeuchts-lemma-in-set-theory}{jdh.hamkins.org/ehrenfeuchts-lemma-in-set-theory}.}}

%
%
%

\maketitle

\begin{abstract}
Ehrenfeucht's lemma~\cite{Ehrenfeucht1973:Discernible} asserts that whenever one element of a model of Peano arithmetic is definable from another, then they satisfy different types. We consider here the analogue of Ehrenfeucht's lemma for models of set theory. The original argument applies directly to the ordinal-definable elements of any model of set theory, and in particular, Ehrenfeucht's lemma holds fully for models of set theory satisfying $V=\HOD$. We show that the lemma can fail, however, in models of set theory with $V\neq\HOD$, and it necessarily fails in the forcing extension to add a generic Cohen real. We go on to formulate a scheme of natural parametric generalizations of Ehrenfeucht's lemma, namely, the principles of the form $\Ehrenfeucht(A,P,Q)$, which asserts that whenever an object $b$ is definable from some $a\in A$ using parameters in $P$, with $b\neq a$, then the types of $a$ and $b$ over $Q$ are different. We also consider various analogues of Ehrenfeucht's lemma obtained by using algebraicity in place of definability, where a set $b$ is \emph{algebraic} in $a$ if it is a member of a finite set definable from $a$ (as in~\cite{HamkinsLeahy:AlgAndImp}). Ehrenfeucht's lemma holds for the ordinal-algebraic sets, we prove, if and only if the ordinal-algebraic and ordinal-definable sets coincide. Using similar analysis, we answer two open questions posed in~\cite{HamkinsLeahy:AlgAndImp}, by showing that (i) algebraicity and definability need not coincide in models of set theory and (ii) the internal and external notions of being ordinal algebraic need not coincide.
\end{abstract}

\section{Introduction}

Ehrenfeucht's lemma asserts that if an element $b$ in a model $M$ of Peano arithmetic ($\PA$) is definable from another distinct element $a$, then the types of $a$ and $b$ in $M$ must be different. This lemma first appeared in a short paper of Ehrenfeucht's \cite{Ehrenfeucht1973:Discernible}, who used it to argue that if a model of $\PA$ is the Skolem closure of a single element $a$, then $a$ is the only element of its type; consequently, such models have no non-trivial automorphisms. Since that time, Ehrenfeucht's lemma has become ubiquitous in the study of models of $\PA$. In light of the extensive transfer of model-theoretic techniques from models of $\PA$ to models of set theory, we find it extremely natural to inquire whether Ehrenfeucht's lemma holds for models of set theory.

The initial answer is that Ehrenfeucht's original argument succeeds directly in models of $V=\HOD$, and more generally, it applies to the ordinal-definable elements of any model of set theory. Nevertheless, we prove that the lemma does not hold in all models of \ZFC, and it definitely fails in the forcing extension obtained by adding a generic Cohen real.

\goodbreak
\begin{thm}\
 \begin{enumerate}
  \item If $a$ is ordinal-definable in a model of set theory $M\models\ZF$ and $b$ is definable from $a$, with $b\neq a$, then $a$ and $b$ satisfy different types in $M$. Consequently, Ehrenfeucht's lemma holds fully in models of $\V=\HOD$.
  \item Meanwhile, Ehrenfeucht's lemma does not hold in all models of set theory. Specifically, if $M$ is any model of $\ZFC$, then in the forcing extension $M[c]$ to add a Cohen real $c$, there are inter-definable elements $a\neq b$ with exactly the same type in $M[c]$. Indeed, there are such elements $a$ and $b$ with $\type^{M[c]}(a/M)=\type^{M[c]}(b/M)$, meaning that $a$ and $b$ satisfy all the same formulas even with parameters from the ground model $M$.
 \end{enumerate}
\end{thm}

These two claims will be proved as theorems~\ref{thm:EhrenfeuchtsLemmaTrueForOD} and~\ref{thm:InterdefinableSetsOfTheSameType}, respectively. Generalizing these observations, we shall introduce in section~\ref{sec:ParametricVersions} the parametric family of principles $\Ehrenfeucht(A,P,Q)$, which holds for a model of set theory $M$ if whenever $a\in A$ and $b$ is definable in $M$ from $a$ using parameters in $P$, with $b\neq a$, then the types of $a$ and $b$ over $Q$ in $M$ are different. So in short, $\Ehrenfeucht(A,P,Q)$ can be expressed by the slogan that \emph{$P$-definability from $A$ implies $Q$-discernibility}. It might be of interest in general to also specify a set $B$ that $b$ has to belong to, but for our present purposes, this is not needed. Usually, the crucial case is that $a$ and $b$ come from the same collection, because if they don't, then they can be distinguished for that very reason, at least for the classes $A$ that will be of interest to us.

These principles unify several natural variations of Ehrenfeucht's lemma that one finds in set theory. The principle $\Ehrenfeucht(M,\emptyset,\emptyset)$, for instance, expresses the original Ehrenfeucht's lemma itself, and we have just mentioned in the theorem above that the principle $\Ehrenfeucht(\OD^M,\emptyset,\emptyset)$ holds in every model of set theory, while  the principle $\Ehrenfeucht(M,\emptyset,\On^M)$ fails in any model $M$ obtained by forcing to add a Cohen real. The principle $\Ehrenfeucht(M,M,\On^M)$ expresses in a model $M$ of set theory that any two distinct elements of $M$ can be distinguished by some formula with an ordinal parameter, which is the Leibniz-Mycielski axiom (see~\cite{Enayat2004:LM}). The principle $\Ehrenfeucht(M,M,\leer)$ holds for a model $M$ if any two distinct elements of $M$ have distinct types, which is precisely the property of $M$ being \emph{Leibnizian} (see~\cite{Enayat2004:Leibnizian}).

Lastly, we shall explore in section~\ref{sec:AlgebraicityVersusDefinability} the relationship between definability and algebraicity, and variations of Ehrenfeucht's lemma which arise by using algebraicity in place of definability.
Specifically, as in~\cite{HamkinsLeahy:AlgAndImp}, a set is algebraic in a parameter if it belongs to a finite set definable from that parameter. We show that Ehrenfeucht's lemma holds on all ordinal-algebraic sets if and only if the ordinal-algebraic sets and the ordinal-definable sets coincide (theorem~\ref{thm:ELonOAequivToOA=OD}), and we also settle several open questions that were asked in~\cite{HamkinsLeahy:AlgAndImp} by pointing out that there are models of set theory in which there are algebraic sets that are not ordinal definable (theorem~\ref{thm:algebraicnondefinable}) and there are models of set theory in which there are objects that are internally but not externally algebraic (theorem~\ref{thm:InternallyButNotExternallyAlgebraic}). Finally, we shall also investigate the algebraic versions of Ehrenfeucht's lemma, the principles $\ExternalAlgebraicEhrenfeucht(A,P,Q)$ and $\InternalAlgebraicEhrenfeucht(A,P,Q)$, which state that if $b$ is (externally/internally) algebraic in $a$ using parameters from $P$, with $b\neq a$, then the types of $a$ and $b$ in $M$ over $Q$ are different -- in short, \emph{$P$-algebraicity from $A$ implies $Q$-discernibility}.

\section{The classic Ehrenfeucht lemma}
\label{sec:classic}

Let us begin by proving Ehrenfeucht's lemma for the ordinal-definable sets in any model of set theory, using essentially Ehrenfeucht's original argument.

\begin{thm}
\label{thm:EhrenfeuchtsLemmaTrueForOD}
If $b$ is definable from an ordinal-definable element $a$ in a model of set theory $M\models\ZF$, with $b\neq a$, then $a$ and $b$ have different types in $M$. $$\type^M(a)\neq\type^M(b)$$
\end{thm}

\begin{proof}
Suppose that $b$ is definable in $M$ from an ordinal-definable element $a$ in $M$. It follows that $b$ also is ordinal-definable, and there is a definable class function $f:M\to M$ such that $f(a)=b$. (Namely, the function defined by $f(x)=y$, if $y$ is the unique set such that $\phi(x,y)$, where $\phi(a,b)$ is a fixed formula defining $b$ from $a$, and otherwise $y=\emptyset$ if there is no such set.) Let us assume that $a<b$ in the definable well-ordering of $\OD^M$; the case $b<a$ can be handled by an essentially identical argument.\footnote{Namely, when $b<a$ one should define the graph to include edges only when $f(x)<x$ and then argue similarly. Alternatively, the case $b<a$ can also be handled more easily as in theorem \ref{thm:LowerRank} by simply counting the number of times one can iterate $f$ before the $\OD$-order stops dropping, and noting that this will be even for $a$ just in case it is odd for $b$.} Let $G$ be the graph on $M$ having as edges the pairs $\{x,f(x)\}$ whenever $x,f(x)\in\OD^M$ and $x<f(x)$ with respect to the \OD\ order. In particular, there is an edge between $a$ and $b$ in $G$. Note that $G$ is loop-free, because it is graded by the $\OD$-order and every node has upward degree at most one, since the only node above $x$ that could be connected directly to $x$ is $f(x)$. So $G$ is a tree, and within each connected component, we have the graph metric $d(x,y)$, which is the length of the shortest path of edges in $G$ connecting $x$ to $y$.
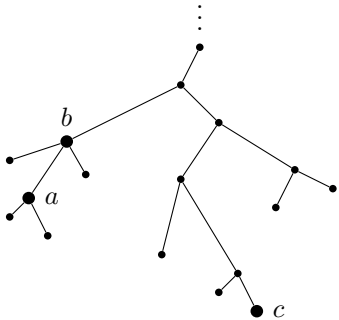
$$\begin{tikzpicture}[scale=.25,yscale=.5]
 \node at (2,10) (a) [circle,fill=black,scale=.5,label=right:$a$] {};
 \node at (4,16) (b) [circle,fill=black,scale=.5,label=above:$b$] {};
 \node at (3,6) (c) [circle,fill=black,scale=.3] {};
 \node at (1,8) (d) [circle,fill=black,scale=.3] {};
 \node at (1,14) (e) [circle,fill=black,scale=.3] {};
 \node at (9,4) (f) [circle,fill=black,scale=.3] {};
 \node at (10,12) (g) [circle,fill=black,scale=.3] {};
 \node at (12,18) (h) [circle,fill=black,scale=.3] {};
 \node at (10,22) (i) [circle,fill=black,scale=.3] {};
 \node at (11,26) (j) [circle,fill=black,scale=.3] {};
 \node at (14,-2) (k) [circle,fill=black,scale=.5,label=right:$c$] {};   
 \node at (13,2) (l) [circle,fill=black,scale=.3] {};
 \node at (15,9) (m) [circle,fill=black,scale=.3] {};
 \node at (16,13) (n) [circle,fill=black,scale=.3] {};
 \node at (18,11) (o) [circle,fill=black,scale=.3] {};
 \node at (12,0) (p) [circle,fill=black,scale=.3] {};
 \node at (5,12.5) (q) [circle,fill=black,scale=.3] {};
 \node at (11,30) (dots) {$\vdots$};
 \draw (c) --(a) --(b) --(i) --(j);
 \draw (d) --(a);
 \draw (e) --(b);
 \draw (k) --(l) --(g) --(h) --(i);
 \draw (f) --(g);
 \draw (m) --(n) --(h);
 \draw (o) --(n);
 \draw (p) --(l);
 \draw (q) --(b);
\end{tikzpicture}$$
Let $c$ be the $\OD$-least element in the connected component of $a$ and $b$ in $G$. Note that $c$ is definable from either $a$ or $b$ by this property. Since $a$ and $b$ are connected by an edge, it follows that the distances $d(a,c)$ and $d(b,c)$ differ by precisely one. In particular, the distance of $a$ to the least member of its connected component is even, just in case the distance of $b$ to the least member of its connected component is odd, and this is therefore a property that distinguishes the types of $a$ and $b$ in $M$. \end{proof}

\begin{cor}
Ehrenfeucht's lemma holds in every model of $\ZF+\V=\HOD$.
\end{cor}

A similar argument applies whenever we define a lower-rank set from a higher-rank set:

\begin{thm}\label{thm:LowerRank}
 If a set $b$ is definable from a set $a$ in a model of set theory $M\models\ZF$ and $b$ has strictly lower rank than $a$, then $a$ and $b$ have different types in $M$.
 $$\type^M(a)\neq\type^M(b)$$
\end{thm}

\begin{proof}
Assume that $b$ is definable from $a$ and has strictly lower rank than $a$ in $M$. As before, there is a definable function $f:M\to M$ such that $b=f(a)$. Starting from $a$, let us iteratively apply the function---producing $a$, $f(a)$, $f(f(a))$, and so on---continuing so long as the ordinal rank of the result continues to drop at each step. Since there is no infinite descending sequence of ordinals in $M$, this process must stop after a finite number of steps in $M$. The number of steps required must be either even or odd in $M$, and since $b=f(a)$ starts one step later than $a$, the number will be even for $a$ if and only if it is odd for $b$, and so the types of $a$ and $b$ must be different.
\end{proof}

Note that in the previous theorem, the rank function could be replaced by any other $M$-definable function. So if $b$ is definable from $a$ in a model $M$ of $\ZF$, and $M\models F(b)\in F(a)$, for some function $F$ which is $M$-definable without parameters, then $a$ and $b$ are discernible. One example would be the cardinality function in a model of $\ZFC$.

\begin{question}
  Does the conclusion of the previous theorem also go through in the case that $b$ has strictly higher rank than $a$?
\end{question}

Meanwhile, Ehrenfeucht's lemma does fail in some models of set theory, and it necessarily fails after forcing to add a generic Cohen real. Indeed, we shall prove that in such a forcing extension, there are inter-definable sets having the same type, even when one allows arbitrary parameters from the ground model.

\begin{thm}\label{thm:InterdefinableSetsOfTheSameType}
If $M$ is a model of \ZFC\ and $M[c]$ is the forcing extension obtained by adding an $M$-generic Cohen real $c$, then in $M[c]$ there are sets $a\neq b$ such that $a$ is definable from $b$ and $b$ is definable from $a$, defined moreover in each case by the same formula, but the types of $a$ and $b$ in $M[c]$ are the same, even when one allows arbitrary parameters from $M$. In other words, $$\type^{M[c]}(a/M)=\type^{M[c]}(b/M).$$
In particular, if \ZFC\ is consistent, then Ehrenfeucht's lemma fails in some models of \ZFC.
\end{thm}

\begin{proof}
It is easy to see that the forcing to add a Cohen real is isomorphic to the $\omega$-fold finite-support product $\P$ of this forcing, and since this will be more convenient, we shall work instead with $\P$. More precisely, conditions in $\P$ are finite
functions $p:\dom(p)\to {}^{{<}\omega}2$, with $\dom(p)\of\omega$, ordered in the natural way by enlarging the domain and extending on each coordinate. Thus, $p(n)$ specifies finitely many bits of the Cohen real to be added at coordinate $n$. Suppose that $G\of\P$ is $M$-generic and consider the forcing extension $M[G]$. For each $n<\omega$, let $g_n=\bigcup\left\{p(n)\st p\in G\right\}$ be the corresponding Cohen real added on coordinate $n$, and let $C=\left\{g_n\st n<\omega\right\}$, $A=\left\{g_{2n}\st n<\omega\right\}$ and $B=\left\{g_{2n+1}\st n<\omega\right\}$. Finally, let
\[ a=\kla{C,A}\quad\text{and}\quad b=\kla{C,B}.\]
Clearly, in $M[G]$ we have $A=C\ohne B$ and $B=C\ohne A$, and so $a$ is definable from $b$, and $b$ is definable from $a$, using the same defining formula in each case.

We shall now show that $a$ and $b$ satisfy the same formulas in $M[G]$, while allowing arbitrary parameters from $M$. Let $\dot{A}$, $\dot{B}$, and $\dot{C}$ be the natural $\P$-names for $A$, $B$, and $C$, respectively. Suppose that $M[G]\models\phi(a,s)$ for some $s\in M$, and fix a condition $p\in\P$ such that $p\forces\varphi(\kla{\dot C,\dot A}, \check s)$. In $M$, we construct an automorphism $\pi$ of $\P$ by swapping certain coordinates as follows. For each $i$ in the support of $p$, which is a finite set, there is by density a coordinate $m_i\in\omega$, with the opposite parity of $i$, such that $p(i)\subseteq g_{m_i}$ and where the $m_i$ are all distinct. Let $\pi$ be an automorphism of $\P$ that interchanges all the even and odd coordinates of the product, undertaken in such as way so as also to swap each $i$ in the support of $p$ with the preselected coordinate $m_i$. Because the support of $p$ is finite, we can find such an automorphism $\pi$ in $M$. It follows that $H=\pi\image G$ is $M$-generic for $\P$ and $M[G]=M[H]$. Our careful choice of $m_i$ ensures specifically that $p\in H$. And since $\pi$ swaps the even and odd coordinates altogether, we have $\dot C^H=C$, $\dot A^H=B$, and $\dot B^H=A$. Since $p\forces\phi(\kla{\dot C,\dot A},\check s)$ and $p\in H$, it now follows that $M[H]\models\phi(\kla{\dot{C}^H,\dot{A}^H},\check{s}^H)$, which means that $M[H]$ and hence $M[G]$ satisfies $\varphi(\kla{C,B},s)$ and thus $\varphi(b,s)$. So we have shown that $\type^{M[G]}(a/M)=\type^{M[G]}(b/M)$, as desired.
\end{proof}

A similar argument is used in \cite[thm~3.1]{Enayat2004:LM} in order to show that the Leibniz-Mycielski axiom can fail in a model of set theory; the new contribution here is that we also make the counterexample sets $a$ and $b$ inter-definable. One interesting consequence of this argument is that $M[c]$ can have no ordinal-definable nor even $M$-definable choice function on the class of all pairs, for we have no definable way to pick one element of $\{a,b\}$ over the other by means of any property in $M[c]$ using parameters from $M$. In particular, there can be no $M$-definable linear ordering of the universe in $M[c]$. A class forcing version of the argument, adding a Cohen set at each regular cardinal, provides a model in which there is no definable linear order of the universe at all, using any set parameters (see \cite{MO110823Hamkins2012:DoesZFCProveTheUniverseIsLinearlyOrderable?}).

\begin{question}
If Ehrenfeucht's lemma holds in a model of set theory $M$, then does it follow that $M\models\V=\HOD$?
\end{question}

\section{A family of parametric Ehrenfeucht principles}
\label{sec:ParametricVersions}

We would like to introduce a natural family of variations of Ehrenfeucht's lemma by allowing parameters in various ways into the definitions and the types. Let us begin with the observation that the original parameter-free version of Ehrenfeucht's lemma is equivalent to a version allowing parameters, which might seem stronger at first sight.

\begin{obs}
\label{obs:EhrenfeuchtWithLocalParameters}
If Ehrenfeucht's lemma holds in a model of set theory $M\models\ZF$, then the following parametric version of it also holds: if $b\neq a$ is definable from $a$ in $M$ using some parameter $p$, then the types of $a$ and $b$ over $p$ in $M$ are different.
\end{obs}

\begin{proof}
If $b$ is definable in $M$ from $a$ with parameter $p$, then the pair $\kla{b,p}$ is definable in $M$ from the pair $\kla{a,p}$, without parameters. Since Ehrenfeucht's lemma holds in $M$, the types of $\kla{a,p}$ and $\kla{b,p}$ in $M$ must be different, and so there is a formula $\varphi$ for which $\varphi(\kla{a,p})$ and $\varphi(\kla{b,p})$ have different truth values in $M$. From this it follows that the types of $a$ and $b$ over $p$ in $M$ are different.
\end{proof}

More interesting generalizations arise when one does not insist that the parameter used to distinguish $a$ and $b$ is the same as the parameter used to define $b$ from $a$.

\begin{defn}
Let $M\models\ZF$, and let $A,P,Q\sub M$. The principle $\Ehrenfeucht(A,P,Q)$ for $M$ asserts that if $a\in A$, $a\neq b$ and $b$ is definable in $M$ from $a$, using parameters from $P$, then there is a formula $\phi$ and parameters $\vc\in Q$ such that $M\models\phi(a,\vc)$, but $M\models\neg\phi(b,\vc)$. So, in plain words, $\Ehrenfeucht(A,P,Q)$ says that $P$-definability from $A$ implies $Q$-discernibility.
We shall write $\Ehrenfeucht(M)$ for $\Ehrenfeucht(M,\leer,\leer)$, which is the original Ehrenfeucht's lemma assertion.
\end{defn}

Expressed in this new terminology, what theorem~\ref{thm:EhrenfeuchtsLemmaTrueForOD} shows is that $\Ehrenfeucht(\OD^M,\leer,\leer)$ holds for every model of set theory $M\models\ZF$. Similarly, theorem~\ref{thm:InterdefinableSetsOfTheSameType} shows that $\Ehrenfeucht(M[c],\leer,M)$ fails if $M[c]$ is the forcing extension of $M\models\ZFC$ by adding a Cohen real $c$. For the rest of this section, we would like to explain how several other principles that have been considered can be expressed in the form $\Ehrenfeucht(A,P,Q)$.

Notice that the principle $\Ehrenfeucht(A,P,Q)$ gets stronger if $A$ or $P$ are enlarged, or if $Q$ is diminished.

\begin{obs}
Fix a model $M\models\ZF$, and a subset $P\sub M$. Then \[\text{if $M$ satisfies}\  \Ehrenfeucht(M),\ \text{then $M$ satisfies}\ \Ehrenfeucht(M,P,P). \]
More generally, if $A\sub M$ and $A\times P^{{<}\omega}\sub A$, then \[\text{if $M$ satisfies}\ \Ehrenfeucht(A,\leer,\leer),\ \text{then $M$ satisfies}\ \Ehrenfeucht(A,P,P).\]
In fact, in each case, it follows that if $b\neq a\in A$ is definable from $a$ using parameters $\vp\in P$ in $M$, then $a$ and $b$ can be distinguished in $M$ using the same parameters $\vp$.
\end{obs}

\begin{proof} The first part follows immediately from Observation~\ref{obs:EhrenfeuchtWithLocalParameters}. A repetition of the proof shows the second part: If $b\neq a$ is definable in $M$ from $a\in A$, using the parameters $\vp\in P$, then $\kla{b,\vp}$ is definable in $M$ from $\kla{a,\vp}$ without parameters. Since $\kla{a,\vp}\in A$, it follows from $\Ehrenfeucht(M,A,\leer,\leer)$ that the $M$-types of $\kla{a,\vp}$ and $\kla{b,\vp}$ are different, and this implies as before that the $M$-types of $a$ and $b$ over $\vp$ are different. \end{proof}

\begin{thm}
Given a model $M\models\ZF$, the principle $\Ehrenfeucht(M,\On^M,\On^M)$, which says that in $M$, ordinal definability implies ordinal discernibility, is first-order expressible. That is, there is a sentence $\phi$ such that for any $M\models\ZF$, $M\models\phi$ if and only if $M$ satisfies $\Ehrenfeucht(M,\On^M,\On^M)$.
\end{thm}

\begin{proof} The sentence in question expresses the following: for any sets $a\neq b$, if $\phi$ is a code for a formula and there are ordinals $\beta<\alpha$ such that $\phi$ defines $b$ from $a$ in $\V_\alpha$, using the parameter $\beta$, then there is a code $\chi(x,y)$ for a formula and there are ordinals $\gamma$ and $\delta$ with $\gamma<\delta$ such that in $\V_\delta$, $\chi(a,\gamma)$ holds, but $\chi(b,\gamma)$ fails. If $\Ehrenfeucht(M,\On^M,\On^M)$ holds for $M$, then the above sentence holds in $M$, by using L\'{e}vy/Montague reflection and using the codes for the actual formulas that define $b$ from $a$ and that distinguish $a$ from $b$. Conversely, if the sentence holds in $M$, then the key observation is that the codes $\phi$ and $\chi$ are ordinals of $M$, even though they may be nonstandard, and thus may be used as additional parameters in the actual formulas defining $b$ from $a$ and distinguishing them.
\end{proof}

The same argument shows:

\begin{cor}
\label{cor:FirstOrderExpressibleEL}
If $A$ is a definable class of $M\models\ZF$, and $B$ and $C$ are definable classes of $M$ containing $\On^M$, then the principle $\Ehrenfeucht(M,A,B,C)$ is first-order expressible over $M$, using the same parameters used to define $A$, $B$ and $C$, if any.
\end{cor}

Next, let's explain that the Leibniz-Mycielski axiom can be expressed as an instance of the principle $\Ehrenfeucht(A,P,Q)$.

\begin{defn}[\cite{Enayat2004:LM}]
A model of set theory $M$ satisfies the Leibniz-Mycielski axiom \LM{}, if whenever $a\neq b$ are sets in $M$, then there is some ordinal $\alpha$ above the ranks of $a$ and $b$ in $M$ and a formula $\varphi$ for which $V_\alpha^M\models\varphi(a)\wedge\neg\varphi(b)$. In other words, $(Th(\V_\alpha,a)\neq Th(\V_\alpha,b))^M$.
\end{defn}

This property is first-order expressible.

\begin{thm}
The Leibniz-Mycielski axiom for a model $M\models\ZF$ is equivalent to the principle $\Ehrenfeucht(M,M,\On^M)$ for $M$, which says that inequality implies ordinal discernibility.
\end{thm}

\begin{proof}
If the Leibniz-Mycielski axiom holds in $M$, then any two elements $a\neq b$ have different types over the ordinals, precisely because $V_\alpha^M\models\varphi(a)\wedge\neg\varphi(b)$, and so $\Ehrenfeucht(M,M,M,\On^M)$ holds. Conversely, if $\Ehrenfeucht(M,M,M,\On^M)$ holds, then what this amounts to is that distinct elements $a\neq b$ of $M$ can be distinguished by formulas using ordinal parameters. So $M\models\varphi(a,\beta)\wedge\neg\varphi(b,\beta)$ for some ordinal parameter $\beta$. By reflection, there is some ordinal $\gamma$ such that $V_\gamma^M$ already sees this. Let $\delta$ be the G\"odel code for $\kla{\gamma,\beta}$ in $M$, and notice that $\delta$ is definable without parameters in $V_{\delta+1}^M$, and so $\gamma$ and $\beta$ are also definable there. So the fact that $V_\gamma^M\models\varphi(a,\beta)$ and $V_\gamma^M\models\neg\varphi(b,\beta)$ are now expressible in $V_{\delta+1}$ as properties distinguishing $a$ and $b$ without any parameters.
\end{proof}

The argument is similar to that of~\cite[Theorem 4.1.1]{Enayat2004:LM}, due to Solovay, showing that \LM{} is equivalent to $\LM^*$, the principle that says that for some fixed formula $\phi(x)$, any two elements $a$ and $b$ can be distinguished in some $\V_\alpha$ by $\phi$.

Enayat points out in~\cite{Enayat2004:LM} that the Leibniz-Mycielski axiom implies that there is a definable linear ordering of the universe, making use of various selection principles. Let us briefly provide a direct argument here, explaining why \LM\ implies that there is a definable bijection of the universe with ${}^{<\On}2$, which is linearly ordered by the lexical order. It suffices by the class version of the Schr\"oder-Cantor-Bernstein theorem to construct a definable injection of the universe into ${}^{<\On}2$. To do so, associate to every object $a$ the set $T_a$ of pairs $\kla{\varphi,\alpha}$ for which $V_\alpha\models\varphi(a)$, for $\alpha$ up to the next $\Sigma_2$-correct ordinal $\beta$ above the rank of $a$. If $a\neq b$, then under $\LM$ there is some ordinal $\alpha$ and formula $\varphi$ for which $V_\alpha\models\varphi(a)\wedge\neg\varphi(b)$, and this must happen before the next $\Sigma_2$-correct ordinal, which implies that $T_a\neq T_b$. Using G\"odel coding, we may view each pair $\kla{\varphi,\alpha}$ as an ordinal, and thus $T_a$ becomes a set of ordinals. So we have a definable injection of the universe into ${}^{<\On}2$ and hence also a definable bijection, as well as a definable linear ordering of the universe. Note that if there is any definable linear ordering of the universe, then there is one that is also set-like, simply by stratifying the universe by rank and using the linear order on each level in turn. Meanwhile, if \ZFC\ is consistent, then there are models of \ZFC\ having no linear ordering of the universe that is definable from parameters (see \cite{MO110823Hamkins2012:DoesZFCProveTheUniverseIsLinearlyOrderable?}).

Note that if $M$ is a model of $\ZF+V=\HOD$, then any two elements of $M$ can be distinguished by an assertion with ordinal parameters (one of them is the $\alpha^{th}$ element in the \OD\ order and the other is not, for some $\alpha$), and so $V=\HOD$ implies $\Ehrenfeucht(M,M,\On^M)$, or, equivalently, $\LM$. Since $\ZF+\V=\HOD$ is equivalent to the existence of ordinal-definable well-ordering of the universe, we are left with the following natural questions, also asked in~\cite{Enayat2004:LM}.

\begin{question}
Which of the implications listed below are strict, and which can be reversed, for a model $M\models\ZF$?
\begin{diagram}[width=6em,height=1.5em,objectstyle=\rm,textflow]
M\models\V=\HOD\\
 \dTo            \\
\text{The Leibniz-Mycielski axiom in }M,\text{ or equivalently, }\Ehrenfeucht(M,M,\On^M)\\
  \dTo           \\
M\ \text{has a definable linear order}\\
\end{diagram}
\end{question}
\bigskip

\noindent Enayat~\cite[conj.~4.3.2]{Enayat2004:LM} conjectures that neither of these implications reverse, and makes several other interesting conjectures concerning further variations and related principles.

Finally, let us consider the principle $\Ehrenfeucht(M,M,\leer)$ for a model $M$ of $\ZF$, which asserts that any two distinct elements $a\neq b$ of $M$ have different types, or, that inequality implies discernibility. This property of a model is known as being \emph{Leibnizian} (see~\cite{Enayat2004:Leibnizian}), following Leibniz's philosophical view on the identity of indiscernibles, so that every object is uniquely determined by the collection of all its properties. This may be viewed as a weak form of pointwise definability (see~\cite{Hamkins2013:PointwiseDefinableModels}): while pointwise definability means that for every object, there is a formula that defines it, being Leibnizian can be viewed as saying that every element of $M$ is defined by the infinite conjunction of all the formulas in its type, i.e., it is definable in the infinitary language with countable conjunctions and disjunctions. So if $M$ is pointwise definable, then $\Ehrenfeucht(M,M,\leer)$ holds for $M$. This implication cannot be reversed, because every pointwise definable model is countable, but Enayat~\cite{Enayat2004:Leibnizian} has proved that there are uncountable Leibnizian models of set theory, which therefore cannot be pointwise definable.

In summary, we have seen that several natural principles studied in set theory can be expressed as instances of the generalized Ehrenfeucht lemma principle $\Ehrenfeucht(A,P,Q)$.

\section{Algebraicity versus Definability}
\label{sec:AlgebraicityVersusDefinability}

We would like now to investigate variants of Ehrenfeucht's lemma that arise by replacing the notion of definability with that of algebraicity. The idea of using the model-theoretic concept of algebraicity in set theory originates in~\cite{HamkinsLeahy:AlgAndImp}. If $M\models\ZF$ and $a,b\in M$, then $b$ is \emph{algebraic} over $a$ in $M$, if there is a finite subset $P$ of $M$ which is definable in $M$ from $a$, such that $b\in P$. Looking at it closely, it becomes clear that there are two ways of understanding the term ``finite'' here: should $M$ believe that $P$ is finite, or should $P$ be finite, as viewed from the outside? The former interpretation is what we call \emph{internal} algebraicity; the latter is \emph{external} algebraicity, and henceforth we shall intend the external meaning unless stated otherwise.

There are two obvious ways in which one can modify Ehrenfeucht's lemma by replacing definability with algebraicity: firstly, one can ask, in a situation where $a,b\in M\models\ZF$ and $b\neq a$ is algebraic in $a$ (possibly allowing parameters), whether it follows that the $M$-type of $b$ is different from the $M$-type of $a$. Secondly, one can ask simply whether the classic Ehrenfeucht's lemma holds for models of set theory not only for the ordinal-definable sets as we proved in theorem \ref{thm:EhrenfeuchtsLemmaTrueForOD}, but also for the ordinal-algebraic sets, the sets that are algebraic in an ordinal parameter. Let us denote the collection of ordinal-algebraic sets of a model $M$ by $\OA(M)$. This latter approach raises a fundamental question on the concept of algebraicity in set theory, which was asked and left open in~\cite{HamkinsLeahy:AlgAndImp}:

\begin{question}[\cite{HamkinsLeahy:AlgAndImp}]\label{question:OA=OD?}
If $M$ is a model of set theory (\ZF{} or \ZFC), does it follow that $\OA(M)=\OD^M$? In other words, are the ordinal-algebraic sets the same as the ordinal-definable sets?
\end{question}

It is interesting that this question remained unanswered, even though a closely related result seems to suggest that $\OA$ should equal $\OD$. In order to state it, let's denote the collection of hereditarily ordinal algebraic sets of $M$ by $\HOA(M)$.

\begin{thm}[{\cite[Theorem 1]{HamkinsLeahy:AlgAndImp}}]
For any model $M$ of $\ZF$, the hereditarily ordinal-algebraic sets are the same as the hereditarily ordinal-definable sets, or in other words, $$\HOA(M)=\HOD^M.$$
\end{thm}

The proof of this result made use of the hereditary nature of the two classes involved, however, and seems to give no information about the relationship between $\OA$ and $\OD$. It was also left open in \cite{HamkinsLeahy:AlgAndImp} whether the parameter-free versions of the concepts of definability and algebraicity may be different.

\begin{question}[\cite{HamkinsLeahy:AlgAndImp}]\label{question:AlgebraicNonDefinablePossible?}
Is there a model of \ZF{} in which there is an element that is algebraic (without parameters) but not definable (without parameters)?
\end{question}

We shall answer all these questions. Let us begin with an observation that shows that (even external) algebraicity does not imply definability, if parameters are allowed.

\begin{obs}
\label{lem:AlgebraicityDoesNotImplyDefinability}
If $c$ is Cohen-generic over $M\models\ZFC$, then there are sets $a$ and $b$ in $M[c]$ such that $b$ is externally algebraic over $a$ in $M[c]$, but $b$ is not definable from $a$, using ordinals, or any parameters from $M$.
\end{obs}

\begin{proof} Let $A$, $B$ and $C$ be the sets used in the proof of theorem~\ref{thm:InterdefinableSetsOfTheSameType}. Let $a=\{A,B\}$, and let $b=A$, say. Clearly, $P=a$ is a finite set definable from $a$, and $b$ is in it. But the proof of theorem~\ref{thm:InterdefinableSetsOfTheSameType} shows that $\kla{\{A,B\},A}$ and $\kla{\{A,B\},B}$ satisfy the same formulas in $M[c]$ (with parameters from $M$). This shows that neither $A$ nor $B$ can be definable from $\{A,B\}$ in $M[c]$, using ordinals, or any parameters from $M$.  \end{proof}

The question underlying any further investigation of whether $\Ehrenfeucht(\OA(M),\leer,\leer)$ for every model $M\models\ZF$ holds is  the fundamental question whether $\OA(M)=\OD^M$ necessarily holds. The following is a sufficient condition for when algebraicity and definability coincide.
\begin{obs}
\label{lem:DefLOImpliesAlgebraicityEqualsDefinability}
If $M\models\ZF$ has a definable linear ordering of its universe, then external algebraicity and definability coincide for $M$. That is, a set $b$ is externally algebraic over $a$ in $M$ just in case it is definable from $a$ in $M$.
\end{obs}

\noindent The fundamental question is answered by the following theorem. The content of the theorem was known to Groszek and Laver, see~\cite{GroszekLaver1987:OD-conjugates}, even though the concept of algebraicity in a set-theoretic context was not considered at the time, and the result was rediscovered by the first author. In order to state it, let us denote Sacks forcing by $\Sacks$, and let us write $\cdegree{x}$ for the constructibility degree of a real $x$.

\begin{thm}
\label{thm:algebraicnondefinable}
If $M\models\ZF+\V=L$ and $a,b$ are mutually $M$-generic Sacks-reals, then $\cdegree{a}\neq\cdegree{b}$ are the only minimal degrees of constructibility in $M[a,b]$, and furthermore, $\cdegree{a}$ and $\cdegree{b}$ have the same type in $M[a,b]$, even allowing parameters from $M$. In particular, the set $\{\cdegree{a},\cdegree{b}\}$ is definable in $M[a,b]$ without parameters, making both $\cdegree{a}$ and $\cdegree{b}$ algebraic (with no parameters) there, but neither $\cdegree{a}$ nor $\cdegree{b}$ are definable in $M[a,b]$, even allowing parameters from $M$. In particular, $\cdegree{a}$ and $\cdegree{b}$ are (parameter-free) algebraic in $M[a,b]$, but not ordinal definable there.
\end{thm}

\begin{proof} (Sketch). Let's just write $L$ for $M$. It is well-known that a Sacks real is minimal, so the degrees of $a$ and $b$ are minimal. A back-and-forth version of the usual argument showing minimality of Sacks reals can be used to show that if $c$ is a real in $L[a,b]$, and $c$ is neither constructible from $a$ nor from $b$, then both $a$ and $b$ are constructible from $c$. This shows that the degrees of $a$ and $b$ are the only minimal degrees in $L[a,b]$.

It follows from a homogeneity property of Sacks forcing $\Sacks$ that the degrees of $a$ and $b$ have the same type in $L[a,b]$ (using ordinals, or arbitrary parameters from $M$): given Sacks conditions $p_0$ and $p_1$, there is an isomorphism between $\Sacks_{\le p_0}$ and $\Sacks_{\le p_1}$, where $\Sacks_{\le r}$ denotes the restriction of the forcing to conditions that are at least as strong as the Sacks condition $r$.

To prove the desired indiscernibility of the constructibility degrees of $a$ and $b$, let $\dot{a}$ and $\dot{b}$ be canonical names for the reals corresponding to the projection of a $\Sacks\times\Sacks$-generic filter on the first and second coordinates, respectively. Suppose $L[a,b]\models\phi(\cdegree{a},\valpha)$, and let $\kla{p,q}\in G\times H$ force this, i.e., in $L$, $\kla{p,q}\forces\phi(\cdegree{\dot{a}},\valpha)$. Let $f:\Sacks_{\le p}\isomorphism\Sacks_{\le q}$ be an isomorphism. Since $G\times H$ is $\Sacks\times\Sacks$-generic over $L$, and since $\kla{p,q}\in G\times H$, it follows that, letting $\bar{G}=G\cap\Sacks_{\le p}$ and $\bar{H}=H\cap\Sacks_{\le q}$, $\bar{G}\times\bar{H}$ is $\Sacks_{\le p}\times\Sacks_{\le q}$-generic over $L$. Let $\bar{I}=f``\bar{G}$, $\bar{J}=f^{-1}``\bar{H}$. Clearly, $\bar{I}\times\bar{J}$ is $\Sacks_{\le q}\times\Sacks_{\le p}$-generic over $L$. Let $I$, $J$ be the closures of $\bar{I}$ and $\bar{J}$ in $\Sacks$ under weakening. Clearly, $I\times J$ is an $\Sacks\times\Sacks$-generic filter for $L$. Let $a'$ and $b'$ be the Sacks reals corresponding to $I$ and $J$, respectively. Because $f$ is an isomorphism, we have that $L[a,b]=L[a',b']$ and because $f\in L$, we have that $\cdegree{a}=\cdegree{a'}$ and $\cdegree{b}=\cdegree{b'}$.
Now, notice that the canonical coordinate switching automorphism of $\Sacks\times\Sacks$, which maps $\kla{r,s}$ to $\kla{s,r}$, induces a transformation of $\Sacks\times\Sacks$ names, which takes $\dot{a}$ to $\dot{b}$ and $\dot{b}$ to $\dot{a}$.
Since we had in $L$ that $\kla{p,q}\forces\phi(\cdegree{\dot{a}},\valpha)$, it follows by applying the coordinate switching isomorphism that $\kla{q,p}\forces\phi(\cdegree{\dot{b}},\valpha)$. Since $\kla{q,p}\in I\times J$, it follows that $L[a',b']\models\phi(\cdegree{b'},\valpha)$. So, since $L[a',b']=L[a,b]$ and $[b']=[b]$, we have that $\varphi(\cdegree{b},\valpha)$ holds in $L[a,b]$. \end{proof}

So it is possible that there are more ordinal algebraic sets than there are ordinal definable ones, which makes the question whether $\Ehrenfeucht(\OA(M),\leer,\leer)$ holds in a model $M\models\ZF$ considerably more attractive. The methods of the previous theorem enable us to settle another question from~\cite{HamkinsLeahy:AlgAndImp} fairly easily: can there be a model of set theory in which there are sets that are internally algebraic but not externally?

\begin{thm}
\label{thm:InternallyButNotExternallyAlgebraic}
If $\ZFC$ is consistent, then there is a model in which there are sets that are internally algebraic, but not externally algebraic, in fact not even externally ordinal algebraic.
\end{thm}

\begin{proof} Let $M$ be a model of $\ZFC+\V=L$, and let $\seq{a_n}{n<\omega}$ be a sequence of Sacks reals over $M$, generic for the $\omega$-fold product of $\Sacks$, with finite-support. Consider the theory consisting of:
\begin{enumerate}
  \item The \ZFC{} axioms,
  \item the sentence expressing ``The set of minimal constructibility degrees is finite,''
  \item for each natural number $n$, the sentence asserting that there are at least $n$ minimal degrees:
  $$\exists v_0\exists v_1\ldots\exists v_{n-1}(\bigwedge_{i<j<n}(v_i\neq v_j)\land\bigwedge_{i<n} \anf{v_i\ \text{is a minimal degree}})$$
  \item\label{item:Indiscernible}
  for every formula $\psi$ with two free variables, the sentence asserting that if one minimal degree satisfies a formula in the parameter $\alpha$, then every minimal degree satisfies that formula in the parameter $\alpha$:
  $$\forall\alpha\in\On\,\forall x((x\text{ is a minimal degree}\,\land\,\psi(x,\alpha))\To\forall y\quad(y\text{ is a minimal degree}\To\psi(y,\alpha)))$$
\end{enumerate}

Using compactness, it is easy to see that this theory is consistent: Given a finite subtheory, one simply has to see what is the largest $n$ such that a sentence of type 3 occurs in it. The methods of theorem~\ref{thm:algebraicnondefinable} show that  $M[a_0,\ldots,a_{n-1}]$ can serve as a model for that subtheory.
In a model of that theory, internally, each minimal degree will be algebraic, since such a model will think that the set of minimal degrees is finite. Note that no nonempty proper subset of the set of all minimal degrees will be ordinal definable in such a model, by the axioms in~\ref{item:Indiscernible}. So no minimal degree will be externally ordinal algebraic, since the set of all minimal degrees in that model is externally infinite. \end{proof}

Let us now formally extend the family of Ehrenfeucht Principles, introduced in section~\ref{sec:ParametricVersions}, to the context of algebraicity. Since the concepts of internal and external algebraicity may differ, we distinguish between them in the formulation of the principles.

\begin{defn}
Let $M\models\ZF$. The external algebraic Ehrenfeucht's lemma for $M$, the principle $\ExternalAlgebraicEhrenfeucht$, says that if $a\neq b$ are elements of $M$ and $b$ is externally algebraic over $a$ in $M$, then $\type^M(a)\neq\type^M(b)$. The internal algebraic Ehrenfeucht's lemma for $M$, the principle $\InternalAlgebraicEhrenfeucht$, gives the same conclusion when $b$ is internally algebraic over $a$ in $M$. More generally, as before, the principles $\ExternalAlgebraicEhrenfeucht(A,P,Q)$ and $\InternalAlgebraicEhrenfeucht(A,P,Q)$ are the corresponding assertions, where $A$ is the set of elements $a$ for which the principle is supposed to hold; $P$ is the set of additional parameters allowed in the algebraic definition of $b$ from $a$; and $Q$ is the set of parameters allowed in order to distinguish the types of $a$ and $b$.
\end{defn}

So the slogan explaining $\ExternalAlgebraicEhrenfeucht(A,P,Q)$ or $\InternalAlgebraicEhrenfeucht(A,P,Q)$ is that (external or internal) \emph{$P$-algebraicity from $A$ implies $Q$-discernibility.}

Clearly, $\InternalAlgebraicEhrenfeucht(M,A,B,C)\implies\ExternalAlgebraicEhrenfeucht(M,A,B,C)$, since external algebraicity implies internal algebraicity. It will turn out that in many cases of interest, the algebraic versions of Ehrenfeucht's lemma are not stronger than the original form.

\begin{lem}
\label{lem:AlgebraicEhrenfeuchtHoldsOnOD}
For any $M\models\ZF$, the principle $\ExternalAlgebraicEhrenfeucht(\OD^M,\leer,\leer)$ for $M$ holds.
\end{lem}

\begin{proof} Let $a\in\OD^M$, and let $b\neq a$ be externally algebraic over $a$ in $M$. So let $P$ be a finite set definable over $M$ from $a$, with $b\in P$. Since $a$ is $\OD^M$, we may assume that $b$ is also $\OD^M$, or else the formula ``$x\in\OD$'' distinguishes $a$ from $b$ over $M$, and we are done. But if $b$ is $\OD^M$ and externally algebraic in $a$, then for some natural number $n$, $b$ is the $n$-th element of $P\cap\OD^M$, in the enumeration according to the well-ordering of $\OD^M$. Since $n$ is an actual natural number, it is definable, and so, $b$ is definable from $a$, not only algebraic in $a$. But then, it follows from Theorem~\ref{thm:EhrenfeuchtsLemmaTrueForOD} that the types of $a$ and $b$ over $M$ differ. \end{proof}

The following theorem answers, among others, our initial question whether Ehrenfeucht's Lemma holds not only on $\OD^M$, but on $\OA(M)$, for a model $M\models\ZF$. It turns out that this is true iff $\OD^M=\OA(M)$.

\begin{thm}
\label{thm:ELonOAequivToOA=OD}
For any model $M\models\ZF$, the following are equivalent:
\begin{enumerate}
  \item
  \label{item:OA=OD}
  $\ExternalOrdinalAlgebraic(M)=\OD^M$

  \item
  \label{item:AEParameterFree}
  $\ExternalAlgebraicEhrenfeucht(\ExternalOrdinalAlgebraic(M),\leer,\leer)$. That is, \\algebraicity from ordinal algebraic sets implies discernibility.

  \item
  \label{item:AE}
  $\ExternalAlgebraicEhrenfeucht(\ExternalOrdinalAlgebraic(M),\leer,\On^M)$. That is, \\algebraicity from ordinal algebraic sets implies ordinal discernibility.

  \item
    \label{item:ELParameterFree}
    $\Ehrenfeucht(\ExternalOrdinalAlgebraic(M),\leer,\leer)$. That is, \\ definability from ordinal algebraic sets implies discernibility.

  \item
    \label{item:EL}
    $\Ehrenfeucht(\ExternalOrdinalAlgebraic(M),\leer,\On^M)$.
    That is, \\definability from ordinal algebraic sets implies ordinal discernibility.
\end{enumerate}
\end{thm}

\begin{proof} The following implications are very easy to see:

\ref{item:OA=OD}$\implies$\ref{item:AEParameterFree} follows from lemma~\ref{lem:AlgebraicEhrenfeuchtHoldsOnOD}.

\ref{item:AEParameterFree}$\implies$\ref{item:ELParameterFree} and~\ref{item:AE}$\implies$\ref{item:EL} are trivial, because if $b$ is definable from $a$, then it is also algebraic in $a$.

\ref{item:AEParameterFree}$\implies$\ref{item:AE} and~\ref{item:ELParameterFree}$\implies$\ref{item:EL} are clear because if $a$ and $b$ can be distinguished without using parameters, then they can also be distinguished using parameters.

The reader is invited to draw a diagram of the implications proven so far. It will become clear that it suffices to prove the following substantial implication.

\ref{item:EL}$\implies$\ref{item:OA=OD}. To see this, suppose $a$ is externally ordinal-algebraic in $M$. Let $P$ be a finite subset of $M$ which is ordinal definable in $M$, and let $a\in P$. Let $Q$ be the set of all linear orderings of the elements of $P$ - i.e., the set of all permutations of $P$. Clearly, $Q$ is $\OD^M$ and finite. Moreover, each linear order is definable from any other one. Fix one such linear ordering of $P$, call it $l$. Since there are only finitely many such linear orderings, $l$ is externally ordinal-algebraic in $M$. By $\Ehrenfeucht(\ExternalOrdinalAlgebraic(M),\leer,\On^M)$ for $M$, for any linear ordering $l'\neq l$ of $P$, there is a formula $\phi_{l'}(x,y)$ and an ordinal $\beta_{l'}\in\On^M$ such that $M\models\phi_{l'}(l,\beta_{l'})\land\neg\phi_{l'}(l',\beta_{l'})$. So $l$ is the unique element $x$ of $Q$  that satisfies $\bigwedge_{l'\in Q\ohne\{l\}}\phi_{l'}(x,\beta_{l'})$ in $M$. So $l$ is ordinal definable in $M$. But then every element of $P$ is the $n$-th element of $P$ in the enumeration $l$, for some $n$. So every element of $P$ is $\OD^M$. The converse is trivial, i.e., every ordinal definable member of $M$ is clearly externally ordinal algebraic. \end{proof}

Let us see what can be said about the internally algebraic versions. Of course, the internal and external versions of algebraicity only differ, if at all, in an $\omega$-nonstandard model.

\begin{lem}
For any model $M\models\ZFC$, the following are equivalent:
\begin{enumerate}
  \item
  \label{item:OA_int=OD}
  $\InternalOrdinalAlgebraic(M)=\OD^M$

  \item
  \label{item:AE_int}
  $\InternalAlgebraicEhrenfeucht(\InternalOrdinalAlgebraic(M),\leer,\On^M)$. That is, \\internal algebraicity from internally ordinal algebraic sets implies ordinal discernibility.

  \item
    \label{item:EL_onInternallyOA}
    $\Ehrenfeucht(\InternalOrdinalAlgebraic(M),\leer,\On^M)$. That is, \\ definability from internally ordinal algebraic sets implies ordinal discernibility.
\end{enumerate}
\end{lem}

\begin{proof} First, note that~\ref{item:OA_int=OD} implies that $\InternalOrdinalAlgebraic(M)=\ExternalOrdinalAlgebraic(M)=\OD^M$, since $\OD^M\sub\ExternalOrdinalAlgebraic(M)\sub\InternalOrdinalAlgebraic(M)$.

\ref{item:OA_int=OD}$\implies$\ref{item:AE_int}: Suppose $a,b\in M$, $b\neq a$, $a\in\InternalOrdinalAlgebraic(M)$, and $b$ is internally algebraic over $a$ in $M$. By~\ref{item:OA_int=OD}, $a\in\OD^M$. So we may assume that $b$ is also $\OD^M$, as otherwise the formula ``$x\in\OD$'' distinguishes the types of $a$ and $b$. But then of course, $a$ and $b$ can be distinguished, using ordinals, since both are definable from ordinals.

\ref{item:AE_int}$\implies$\ref{item:EL_onInternallyOA} is trivial, since if $b$ is definable from $a$, then $b$ is also internally algebraic in $a$.

\ref{item:EL_onInternallyOA}$\implies$\ref{item:OA_int=OD}:
Suppose $a\in M$ is internally ordinal algebraic. Let $a\in P$, $P$ $M$-finite, $P\in\OD^M$. We have to show that $a$ is ordinal definable in $M$. As before, and working inside $M$, let $Q$ be the set of linear orderings of $P$. Clearly, $Q$ is $M$-finite and $\OD^M$, since $P$ is. Fix one particular linear ordering $l\in Q$ (so $l$ is internally ordinal algebraic). Let us fix the number of elements of $P$, say $n$ (a possibly nonstandard number in $M$). For any $l'\in Q$, different from $l$, there is a permutation $\pi_{l'}$ of the first $n$ natural numbers such that $l'$ is the result of applying $\pi_{l'}$ to $l$. Of course, $\pi_{l'}$ can be viewed as a (nonstandard) number of $M$. Now, note that the pair $\kla{\pi_{l'},l'}$ is definable from the pair $\kla{\pi_{l'},l}$. So by $\Ehrenfeucht(\InternalOrdinalAlgebraic(M),\leer,\On^M)$ for $M$, there is a formula $\phi_{l'}$ such that in $M$, $\phi_{l'}(\pi_{l'},l)$ holds but $\phi_{l'}(\pi_{l'},l')$ fails. By Levy reflection, this can be seen locally in $M$, and it follows that for every $l'\in Q$, $l'\neq l$, there is an $M$-ordinal $\alpha_{l'}$ and an $M$-formula $\psi_{l'}$ such that in $M$, for every $l'\in Q$ with $l'\neq l$, $\V_{\alpha_{l'}}\models\psi_{l'}(\pi_{l'},l)\land\neg\psi_{l'}(\pi_{l'},l')$.
In fact, the function $l'\mapsto\kla{\pi_{l'},\psi_{l'},\alpha_{l'}}$ can be arranged to belong to $M$. In $M$, let us enumerate $\{\kla{\pi_{l'},\psi_{l'},\alpha_{l'}}\st l'\in Q\ohne\{l\}\}$ as $\{\kla{\pi_r,\psi_r,\alpha_r}\st r<p\}$, for some (possibly nonstandard) natural number $p$. Consider the $M$-formula
\[\bigwedge_{i<p}\anf{\V_{\alpha_i}\models\psi_{i}(\pi_i,x)}\]
in the parameter $c:=\kla{\kla{\pi_0,\psi_0,\alpha_0},\ldots,\kla{\pi_{p-1},\psi_{p-1},\alpha_{p-1}}}$.
Since the whole parameter $c$ can be coded by a single ordinal $\theta$ in $M$, that formula can be expressed as $\chi(\theta,x)$, for an actual external formula $\chi$. It follows then that $l$ is the unique member of $Q$ such that $M\models\chi(\theta,l)$. Since $Q$ was $\OD^M$, it follows that $l$ is $\OD^M$. But then, $a$ is the $s$-th element of $Q$ according to the enumeration of $Q$ by $l$, for some (possibly nonstandard) number $s$. So $a$ is ordinal definable in $M$. \end{proof}

\begin{question}
If $\OD^M=\InternalOrdinalAlgebraic(M)$, does it follow that $\InternalAlgebraicEhrenfeucht(\InternalOrdinalAlgebraic(M),\leer,\leer)$ holds for $M$?
\end{question}

\begin{defn}
  Let $M\models\ZF$. We define that $\Definable(M)$ is the collection of $a\in M$ that are definable over $M$ without parameters and that $\ExternalAlgebraic(M)$ is the collection of $a\in M$ that are externally algebraic in $M$, without parameters.  The collection $\InternalAlgebraic(M)$ is defined similarly.
\end{defn}

The same proof shows the following parameter-free version of the previous theorem.

\begin{thm}
For any model $M\models\ZF$, the following are equivalent:
\begin{enumerate}
  \item
  \label{item:A=D}
  $\ExternalAlgebraic(M)=\Definable(M)$, that is, algebraicity is the same as definability.
  \item
  \label{item:AEA}
  $\ExternalAlgebraicEhrenfeucht(\ExternalAlgebraic(M),\leer,\leer)$. That is, algebraicity in an algebraic set implies discernibility.
  \item
    \label{item:ELA}
    $\Ehrenfeucht(\ExternalAlgebraic(M),\leer,\leer)$. That is, definability from an algebraic set implies discernibility.
\end{enumerate}
\end{thm}

\end{document}